\documentstyle[twoside,11pt]{article}
  \newcommand{\bref}[5]{\noindent\parbox[t]{0.8cm}{#1}\parbox[t]{15.2cm}{{#2}{\it #3}{\bf #4}#5}\par}

  \def\sg{\sigma}

\begin{document}
  \title{\bf {A P\'{o}sa-type condition of potentially $_3C_\ell$-graphic sequences}
  \thanks{Supported by National Natural Science Foundation of China (No. 11961019).}}
  \author{Guang-Ming Li, Jian-Hua Yin\thanks{Corresponding author.}\\
   {\small School of Science, Hainan University, Haikou 570228, P.R. China}\\
{\small E-mail: yinjh@hainanu.edu.cn}}
\date{}
\maketitle

\begin{center}
\parbox{0.9\hsize}
{\small {\bf Abstract.}\ \ A non-increasing sequence
$\pi=(d_1,\ldots,d_n)$ of nonnegative integers is said to be {\it
graphic} if it is realizable by a simple graph $G$ on $n$ vertices.
A graphic sequence $\pi=(d_1,\ldots,d_n)$ is said to be {\it
potentially} $_3C_\ell$-{\it graphic} if there is a realization of
$\pi$ containing cycles of every length $r$, $3\le r\le \ell$. It is
well-known that if the non-increasing degree sequence
$(d_1,\ldots,d_\ell)$ of a graph $G$ on $\ell$ vertices satisfies
the P\'{o}sa condition that $d_{\ell+1-i}\ge i+1$ for every $i$ with
$1\le i<\frac{\ell}{2}$, then $G$ is either pancyclic or bipartite.
In this paper, we obtain a P\'{o}sa-type condition of potentially
$_3C_\ell$-graphic sequences, that is, we prove that if $\ell\ge 5$
is an integer, $n\ge \ell$ and $\pi=(d_1,\ldots,d_n)$ is a graphic
sequence with $d_{\ell+1-i}\ge i+1$ for every $i$ with $1\le
i<\frac{\ell}{2}$, then $\pi$ is potentially $_3C_\ell$-graphic.
This result improves a Dirac-type condition of potentially
$_3C_\ell$-graphic sequences due to Yin et al. [Appl. Math. Comput.,
353 (2019) 88--94], and asymptotically answers a problem due to Li
et al. [Adv. Math., 33 (2004) 273--283]. As an application, this
result also completely implies the value $\sg(C_\ell,n)$ for
$\ell\ge 5$ and $n\ge \ell$, improving the result of Lai [J. Combin.
Math. Combin. Comput., 49 (2004) 57--64].\\

}
\end{center}

\section*{1. Introduction}
\hskip\parindent Graphs in this paper are finite, undirected and
simple. Terms and notation not defined here are from [1]. A sequence
$\pi=(d_1,\ldots,d_n)$ of nonnegative integers (not necessarily in
non-increasing order) is said to be {\it graphic} if it is the
degree sequence of a simple graph $G$ on $n$ vertices. In this case,
$G$ is referred to as a {\it realization} of $\pi$. The set of all
sequences $(d_1,\ldots,d_n)$ of nonnegative integers with $d_1\ge
\cdots \ge d_n$ and $d_1\le n-1$ is denoted by $NS_n$. The set of
all graphic sequences in $NS_n$ is denoted by $GS_n$. For a
nonnegative integer sequence $\pi=(d_1,\ldots,d_n)$, we denote
$\sigma(\pi)=d_1+\cdots+d_n$. A sequence $\pi\in GS_n$ is said to be
{\it potentially} $C_\ell$-{\it graphic} if there is a realization
of $\pi$ containing $C_\ell$ as a subgraph, where $C_\ell$ is the
cycle of length $\ell$. Moreover, a sequence $\pi\in GS_n$ is said
to be {\it potentially} $_3C_\ell$-{\it graphic} if there is a
realization of $\pi$ containing cycles of every length $r$, $3\le
r\le \ell$. Li and Yin proposed two problems as follows.

{\bf Problem 1.1} [8]\ \ {\it Give a criteria of potentially
$C_\ell$-graphic sequences for $\ell\ge 3$.}

{\bf Problem 1.2} [8]\ \ {\it Give a criteria of potentially
$_3C_\ell$-graphic sequences for $\ell\ge 4$.}

In [11], Yin gave a solution to Problem 1.1. However, Problem 1.2
seems to be much more challenging over time. Chen et al.
investigated the condition of potentially $_3C_\ell$-graphic
sequences for $\ell=5,6$. Note that the symbol $x^y$ in a sequence
stands for $y$ consecutive terms, each equal to $x$.

{\bf Theorem 1.1} [2,3]\ \ {\it Let $\pi=(d_1\ldots,d_n)\in GS_n$.

(i)\ \ If $n\ge 5$, $d_5\ge 2$ and $d_4\ge 3$, then $\pi$ is
potentially $_3C_5$-graphic.

(ii)\ \ If $n\ge 6$, $d_6\ge 2$ and $d_5\ge 3$, then $\pi$ is
potentially $_3C_6$-graphic.}

Recently, Yin et al. established a Dirac-type condition of
potentially $_3C_\ell$-graphic sequences for $\ell\ge 5$ as follows.

{\bf Theorem 1.2} [12]\ \ {\it If $\ell\ge 5$ is an integer, $n\ge
\ell$ and $\pi=(d_1,\ldots, d_n)\in GS_n$ with $d_\ell\ge
\lceil\frac{\ell}{2}\rceil$, then $\pi$ is potentially
$_3C_\ell$-graphic.}

It is well-known that if the non-increasing degree sequence
$(d_1,\ldots,d_\ell)$ of a graph $G$ on $\ell$ vertices satisfies
the P\'{o}sa condition that $d_{\ell+1-i}\ge i+1$ for every $i$ with
$1\le i<\frac{\ell}{2}$, then $G$ is either pancyclic or bipartite
(see [10]). Motivated by this result, we obtain a P\'{o}sa-type
condition of potentially $_3C_\ell$-graphic sequences for $\ell\ge
5$.

{\bf Theorem 1.3}\ \ {\it If $\ell\ge 5$ is an integer, $n\ge \ell$
and $\pi=(d_1,\ldots,d_n)\in GS_n$ with $d_{\ell+1-i}\ge i+1$ for
every $i$ with $1\le i<\frac{\ell}{2}$ (i.e.,
$i=1,\ldots,\lceil\frac{\ell}{2}\rceil-1$), then $\pi$ is
potentially $_3C_\ell$-graphic.}

Clearly, Theorem 1.3 improves Theorem 1.2. Moreover, Theorem 1.3 is
best possible in the sense that the condition $d_{\ell+1-i}\ge i+1$
cannot be replaced by $d_{\ell+1-i}\ge i$ for every
$i=1,\ldots,\lceil\frac{\ell}{2}\rceil-1$, since both
$((n-1)^i,(2m-i-1)^{m-i},m^{m-i},i^{n-2m+i})$ for $\ell=2m$ and
$((n-1)^i,(2m-i)^{m-i+1},(m+1)^{m-i},i^{n-2m+i-1})$ for $\ell=2m+1$
are graphic, but $((n-1)^i,(2m-i-1)^{m-i},m^{m-i},i^{n-2m+i})$ is
not potentially $_3C_{2m}$-graphic and
$((n-1)^i,(2m-i)^{m-i+1},(m+1)^{m-i},i^{n-2m+i-1})$ is not
potentially $_3C_{2m+1}$-graphic. This implies that Theorem 1.3
answers Problem 1.2 asymptotically.

In [7], Lai considered an extremal problem of potentially
$C_\ell$-graphic sequences as follows. Determine the minimum even
integer $\sg (C_\ell,n)$ such that every sequence $\pi\in GS_n$ with
$\sg (\pi)\ge \sg (C_\ell,n)$ is potentially $C_\ell$-graphic. We
will refer to $\sg (C_\ell,n)$ as the {\it potential number} of
$C_\ell$. As $\sg(\pi)$ is twice the number of edges in any
realization of $\pi$, this problem can be viewed as a potential
degree sequence relaxation of the Tur\'{a}n number of $C_\ell$. Lai
[7] determined $\sg(C_{2m+1},n)=m(2n-m-1)+2$ for $m\ge 2$ and $n\ge
3m$, and $\sg(C_{2m+2},n)=m(2n-m-1)+4$ for $m\ge 2$ and $n\ge 5m-2$.
As an application of Theorem 1.3, we will completely determine the
value $\sg(C_\ell,n)$ for $\ell\ge 5$ and $n\ge \ell$, improving the
above result of Lai.

{\bf Theorem 1.4}\ \ Let $m\ge 2$. Then
$\sg(C_{2m+1},n)=\max\{2n+4m^2-6m,m(2n-m-1)\}+2$ for $n\ge 2m+1$ and
$\sg(C_{2m+2},n)=\max\{2n+4m^2-2m-2,m(2n-m-1)+2\}+2$ for $n\ge
2m+2$.

\section*{2. Proofs of Theorem 1.3--1.4}
\hskip\parindent We first introduce some useful known theorems and
lemmas.

Let $\pi=(d_1,\ldots,d_n)\in NS_n$, and let $k$ be an integer with
$1\le k\le n$. We define $\pi_k'$ to be a non-increasing sequence so
that $\pi_k'$ is obtained from $\pi$ by deleting $d_k$ and
decreasing $d_k$ largest terms from
$d_1,\ldots,d_{k-1},d_{k+1},\ldots,d_n$ each by one unity. We say
that $\pi_k'$ is the {\it residual sequence} obtained from $\pi$ by
laying off $d_k$. It is easy to see that if $\pi_k'$ is graphic then
so is $\pi$, since a realization of $\pi$ can be obtained from a
realization of $\pi_k'$ by adding a new vertex that is adjacent to
those vertices whose degrees are reduced by one unity in going from
$\pi$ to $\pi_k'$. In fact, more is true:

{\bf Theorem 2.1} [6]\ \ {\it $\pi\in GS_n$ if and only if
$\pi_k'\in GS_{n-1}$.}

{\bf Theorem 2.2} [4]\ \ {\it Let $\pi=(d_1,\ldots,d_n)\in NS_n$,
where $\sigma(\pi)$ is even. Then $\pi \in GS_n$ if and only if
$\sum\limits_{i=1}^{t}d_i\le
t(t-1)+\sum\limits_{i=t+1}^{n}\min\{t,d_i\}$ for each $t$ with $1\le
t\le n-1$.}

{\bf Remark.}\ \ Nash-Williams [9] proved that to show that
$\pi=(d_1,\ldots,d_n)\in NS_n$ is graphic, it is enough to check
$\sum\limits_{i=1}^td_i \le t(t-1)+\sum\limits_{i=t+1}^n
\min\{t,d_i\}$ for those $t$ which satisfy $d_t>d_{t+1}$. Moreover,
it is not difficult to prove that to show that
$\pi=(d_1,\ldots,d_n)\in NS_n$ is graphic, it is enough to check
$\sum\limits_{i=1}^td_i \le t(t-1)+\sum\limits_{i=t+1}^n
\min\{t,d_i\}$ for each $t$ with $1\le t\le f(\pi)$, where
$f(\pi)=\max\{i|d_i\ge i\}$.

{\bf Theorem 2.3} [5]\ \ {\it If $\pi=(d_1,\ldots,d_n)\in GS_n$ has
a realization $G$ containing $H$ as a subgraph, then there exists a
realization $G'$ of $\pi$ containing $H$ on those vertices with
degrees $d_1,\ldots,d_{|V(H)|}$.}

{\bf Lemma 2.1} [12]\ \ {\it If $\ell\ge 3$ and
$\pi=(d_1,\ldots,d_\ell)\in GS_\ell$ with $d_{\ell+1-i}\ge i+1$ for
every $i$ with $1\le i<\frac{\ell}{2}$ (i.e.,
$i=1,\ldots,\lceil\frac{\ell}{2}\rceil-1$), then $\pi$ has a
pancyclic realization, that is, $\pi$ has a realization containing
cycles of all possible length from $3$ up to $\ell$.}

{\bf Lemma 2.2} [12]\ \ {\it If $m\ge 4$, $n\ge 2m$ and
$\pi=(d_1,\ldots, d_n)\in GS_n$ with $d_{2m}\ge m$, then $\pi$ has a
realization containing $C_3,\ldots, C_{2m}$ on those vertices with
degrees $d_1,\ldots,d_{2m}$.}

We also need a series of lemmas.

{\bf Lemma 2.3}\ \ {\it Let $m\geq 5$, $n\geq 2m+1$ and
$\pi=(d_1,\ldots,d_n)\in GS_n$ with $m-1\ge d_n\geq 1$,
$d_{d_n-1}\geq m+1$, $d_{d_n}=\cdots=d_{m+1}=m$ and $d_{2m+1-i}\geq
i+1$ for $1\leq i\leq m-1$. If $n\geq 2m+2$, let
$\omega=(d_1-1,\ldots,d_{d_n-1}-1,d_{d_n},\ldots,d_{n-2},d_{n-1}-1)$
(that is, $\omega$ is obtained from $\pi$ by deleting the term
$d_n$, and then decreasing the term $d_i$ by one unity for each
$i\in \{1,\ldots,d_n-1,n-1\}$), and if $n=2m+1$ and there exists an
integer $k$ with $1\le k\le m-2$ so that $d_{2m+1-k}\geq k+2$ (we
choose $k$ to be minimum), let
$\omega=(d_1-1,\ldots,d_{d_n-1}-1,d_{d_n},\ldots,d_{2m-k},d_{2m+1-k}-1,d_{2m+2-k},\ldots,d_{n-1})$,
then $\omega$ is graphic.}

{\bf Proof.}\ \ Denote $\omega=(\rho_1,\ldots,\rho_{n-1})$. Clearly,
$\rho_i=d_i-1$ for $1\leq i\leq d_n-1$,
$\rho_{d_n}=\cdots=\rho_{m+2}=m$ and $\rho_{n-1}\geq d_n-1$. For
$1\leq t\leq d_n-1$, we have $\sum\limits_{i=1}^{t}\rho_i\leq
t(n-2)=t(t-1) + t(n-1-t)\leq t(t-1) +
\sum\limits_{i=t+1}^{n-1}\min\{t,\rho_i\} $. Moreover, for $t=m+2$,
we also have $\sum\limits_{i=1}^t\rho_i\leq
(n-2)(d_n-1)+m(m+3-d_n)\leq (m+2)(m+1) + (d_n-1)(n-t-1) \leq t(t-1)
+ \sum\limits_{i=t+1}^{n-1}\min\{t,\rho_i\}$. Thus by
$f(\rho)=m<m+2$ and Theorem 2.2, $\omega$ is graphic.\ \ $\Box$

{\bf Lemma 2.4}\ \ {\it Let $m\ge 5$, $n\ge 2m+1$ and
$\pi=(d_1,\ldots,d_n)\in GS_n$ with $m-1\ge d_n\geq 1$,
$d_{d_n-1}=d_{d_n}=\cdots= d_{m+2}=m$ and
$d_{m+3}=\cdots=d_{2m+1}=m-1$. Denote $p=\max\{i|d_i\geq m+1\}$. Let
$\rho=(d_1,\ldots,d_{m+3},d_{m+4}-1,\ldots,d_{2m}-1,d_{2m+2}-1,d_{2m+3}-1,d_{2m+4},\ldots,d_n)$
when $p=0$ and $n\geq 2m+3$ (that is, $\rho$ is obtained from $\pi$
by deleting the term $d_{2m+1}=m-1$, and then decreasing the term
$d_i$ by one unity for each $i\in \{m+4,\ldots,2m,2m+2,2m+3\}$),
$\rho=(d_1-1,d_2,\ldots,d_{m+3},d_{m+4}-1,\ldots,d_{2m}-1,d_{2m+2}-1,d_{2m+3},\ldots,d_n)$
when $p=1$ and $n\geq 2m+2$,
$\rho=(d_1-1,d_2-1,d_3,\ldots,d_{m+3},d_{m+4}-1,\ldots,d_{2m}-1,d_{2m+2},\ldots,d_n)$
when $p=2$, and
$\rho=(d_1-1,\ldots,d_p-1,d_{p+1},\ldots,d_{m+3},d_{m+4}-1,\ldots,d_{2m+2-p}-1,d_{2m+3-p},\ldots,d_{2m},d_{2m+2},\ldots,d_n)$
when $p\ge 3$. Then $\rho$ is graphic.}

{\bf Proof.}\ \ Let $\rho'=(\rho_1,\ldots,\rho_{n-1})$ be a
rearrangement in non-increasing order of the $n-1$ terms of $\rho$.
Clearly, $\rho_i=d_i-1$ for $1\le i\le p$, $\rho_i=d_i=m$ for
$p+1\le i\le m+2$ and $\rho_{n-1}\ge d_n-1$. For $1\le t\le p$, by
$p\le d_n-2$, we have $\sum\limits_{i=1}^t \rho_i\leq t(n-2)=t(t-1)
+ t(n-1-t)=t(t-1)+\sum\limits_{i=t+1}^{n-1}\min\{t,\rho_i\}$.
Moreover, for $t=m+2$, we also have $\sum\limits_{i=1}^t \rho_i\leq
p(n-2)+(t-p)m\le t(t-1) +(n-1-t)p\le t(t-1)+
\sum\limits_{i=t+1}^{n-1}\min\{t,\rho_i\}$. Thus by $f(\rho')=m<m+2$
and Theorem 2.2, $\rho'$ is graphic, that is, $\rho$ is graphic. \ \
$\Box$

{\bf Lemma 2.5}\ \ {\it Let $m\ge 5$, $n\geq 2m+1$ and
$\pi=(d_1,\ldots,d_n)\in GS_n$ with $m-1\ge d_{2m}\ge d_{2m+1}\ge
1$, $d_{d_n-1}=d_{d_n}=\cdots=d_{m+1}=m$ and $d_{2m+1-i}\geq i+1$
for $1\leq i\leq m-1$. Denote $p=\max\{i|d_i\geq m+1\}$. Clearly,
$d_{2m}\ge p+2$. Let
$\omega_1=(d_1-1,\ldots,d_m-1,d_{m+1},d_{m+3},\ldots,d_{n})$ (that
is, $\omega_1$ is obtained from $\pi$ by deleting the term
$d_{m+2}=m$, and then decreasing the term $d_i$ by one unity for
each $i\in \{1,\ldots,m\}$). In fact, $\omega_1$ (not in
non-increasing order) is the residual sequence obtained from $\pi$
by laying off $d_{m+2}=m$. If $p\ge 1$ or $d_{2m}\le m-2$, let
$\omega_2=(d_1-2,\ldots,d_p-2,d_{p+1}-1,\ldots,d_{m+1}-1,d_{m+3}-1,\ldots,d_r-1,d_{r+1},\ldots,d_{2m-1},d_{2m+1},\ldots,d_{n})$,
where $r=m+1+d_{2m}-p$. If $p=0$ and $d_{2m}=m-1$, let
$\omega_2=(d_1-2,\ldots,d_p-2,d_{p+1}-1,\ldots,d_{m+1}-1,d_{m+3}-1,\ldots,d_{2m-1}-1,d_{2m+1}-1,d_{2m+2},\ldots,d_{n})$.
Then $\omega_2$ is graphic.}

{\bf Proof.}\ \ If $d_r=m$, then $\omega_2$ (not necessarily in
non-increasing order) is the residual sequence obtained from
$\omega_1$ by laying off $d_{2m}$. By Theorem 2.1, both $\omega_1$
and $\omega_2$ are graphic. Assume $d_r\le m-1$. Let
$\rho=(\rho_1,\ldots,\rho_{n-2})$ be a rearrangement in
non-increasing order of the $n-2$ terms of $\omega_2$. Clearly,
$\rho_k=d_k-2$ for $1\leq k\leq p$, $\rho_k=m-1$ for $p+1\leq k\leq
m+1$ and $\rho_{n-2}\geq d_n-1$. For $1\leq t\leq p$, by $p\leq
d_n-1$, we have $ \sum\limits_{i=1}^{t}\rho_i=
\sum\limits_{i=1}^{t}(d_i-2)\leq t(n-3)=t(t-1)+t(n-2-t)\leq t(t-1) +
\sum\limits_{i=t+1}^{n-2}\min\{t,\rho_i\}$. Moreover, for $t=m+1$,
we also have $\sum\limits_{i=1}^t \rho_i\leq p(n-3)+(t-p)(m-1)\leq
t(t-1) +p(n-2-t)\leq  t(t-1)+
\sum\limits_{i=t+1}^{n-2}\min\{t,\rho_i\}$. Thus by
$f(\rho)=m-1<m+1$ and Theorem 2.2, $\rho$ is graphic, that is,
$\omega_2$ is graphic.  \ \ $\Box$

{\bf Lemma 2.6} \ \ {\it Let $m\geq 5$ and
$\pi=(d_1,\ldots,d_{2m+1})\in GS_{2m+1}$ with $d_i=m$ for $2\leq
i\leq m+2$ and $d_{2m+1-i}=i+1$ for $1\leq i\leq m-2$. If
$d_{2m+1}=2$ and $d_1\geq m+1$, or if $d_{2m+1}=1$ and $d_1=m$, then
$\pi$ has a realization containing $C_3,\ldots, C_{2m}$ on those
vertices with degrees $d_1,\ldots,d_{2m}$.}

{\bf Proof.}\ \ Firstly, we assume that $d_{2m+1}=2$ and $d_1\geq
m+2$.

Then $\pi=(d_1,m^{m+1},m-1,m-2,\ldots,2,2)$. Let $\omega_0=\pi$. Let
$\omega_1=(d_1',\ldots,d_{2m}')$ be the residual sequence obtained
from $\pi$ by laying off the term $d_{2m+1}=2$. Then
$\omega_1=(d_1-1,m^m,m-1,m-1,m-2,\ldots,2)$. Let
$\omega_2=(d_1'',\ldots,d_{2m-1}'')$ be the residual sequence
obtained from $\omega_1$ by laying off the term $d_{m+3}'=m-1$. Then
$\omega_2=(d_1-2,m^2,(m-1)^{m-1},m-2,\ldots,2)$. Let
$\omega_3=(d_1''',\ldots,d_{2m-2}''')$ be the residual sequence
obtained from $\omega_2$ by laying off the term $d_{2m-2}''=3$. Then
$\omega_3=(d_1-3,(m-1)^{m+1},m-2,\ldots,2)$. By Theorem 2.1 and
Lemma 2.1, $\omega_3$ has a pancyclic realization $G_3$. For
$i=2,1,0$ in turn, let $G_i$ be the graph obtained from $G_{i+1}$ by
adding a new vertex $u_i$ that is adjacent to those vertices whose
degrees are reduced by one unity in going from $\omega_i$ to
$\omega_{i+1}$. Then $G_0$ is a realization of $\pi$ and contains
$G_3$. Let $C=x_1x_2\ldots x_{2m-2}x_1$ be a Hamilton cycle of
$G_3$. Then $V(G_0)=\{x_1,\ldots,x_{2m-2},u_0,u_1,u_2\}$,
$d_{G_0}(u_0)=d_{2m+1}=2$, $d_{G_0}(u_1)=d_{m+3}=m-1$ and
$d_{G_0}(u_2)=d_{2m-1}=3$. Thus $C$ lies on those vertices with
degrees $d_1,\ldots,d_{m+2},d_{m+4},\ldots,d_{2m-2},d_{2m}$. For
convenience, let $\{i_1,\ldots,i_{m+1}\}\subseteq
\{1,2,\ldots,2m-2\}$ so that $d_{G_0}(x_{i_j})=d_j$ for $1\leq j\leq
m+1$. Clearly, $u_2x_{i_1},u_2x_{i_2},u_2x_{i_3}\in E(G_0)$ and
$u_1x_{i_j}\in E(G_0)$ for $j=1$ and $4\leq j\leq m+1$. Moreover, it
is easy to see that there exist four distinct positive integers
$p_1,p_2,p_3,p_4\in \{i_1,\ldots,i_{m+1}\}$ so that
$x_{p_1}x_{p_2},x_{p_3}x_{p_4}\in E(C)$. We now show that $\pi$ has
a realization containing $C_3,\ldots, C_{2m}$ on those vertices with
degrees $d_1,\ldots,d_{2m}$. Without loss of generality, we only
need to consider two cases.

{\bf Case 1}\ \ $x_{i_1}\neq x_{p_1},x_{p_2},x_{p_3},x_{p_4}$.

If $x_{i_2}=x_{p_1}$ and $x_{i_3}=x_{p_2}$, then
$u_2x_{p_1},u_2x_{p_2},u_1x_{p_3},u_1x_{p_4}\in E(G_0)$, implying
that $G_0$ contains $C_3,\ldots,C_{2m}$ on those vertices with
degrees $d_1,\ldots,d_{2m}$. If $x_{i_2}=x_{p_1}$ and $x_{i_3}\neq
x_{p_2}$, let
$G_1=G_0-\{u_2x_{i_3},u_1x_{p_2}\}+\{u_2x_{p_2},u_1x_{i_3}\}$, then
$G_1$ is a realization of $\pi$ and
$u_2x_{p_1},u_2x_{p_2},u_1x_{p_3},u_1x_{p_4}\in E(G_1)$. Thus $G_1$
contains $C_3,\ldots,C_{2m}$ on those vertices with degrees
$d_1,\ldots,d_{2m}$. If $x_{i_2}\neq x_{p_1},x_{p_2}$ and
$x_{i_3}\neq x_{p_1},x_{p_2}$, let
$G_1=G_0-\{u_2x_{i_2},u_2x_{i_3},u_1x_{p_1},u_1x_{p_2}\}+\{u_2x_{p_1},u_2x_{p_2},u_1x_{i_2},u_1x_{i_3}\}$,
then $G_1$ is a realization of $\pi$ and contains
$C_3,\ldots,C_{2m}$ on those vertices with degree
$d_1,\ldots,d_{2m}$.

{\bf Case 2}\ \ $x_{i_1}=x_{p_1}$.

If $x_{i_2}=x_{p_3}$ and $x_{i_3}=x_{p_4}$, then
$u_2x_{p_3},u_2x_{p_4},u_1x_{p_1},u_1x_{p_2}\in E(G_0)$ and $G_0$
contains $C_3,\ldots,C_{2m}$ on those vertices with degrees
$d_1,\ldots,d_{2m}$. If $x_{i_2}=x_{p_3}$ and $x_{i_3}\neq x_{p_4}$,
let $G_1=G_0-\{u_2x_{i_3},u_1x_{p_4}\}+\{u_2x_{p_4},u_1x_{i_3}\}$,
then $G_1$ is a realization of $\pi$ and contains
$C_3,\ldots,C_{2m}$ on those vertices with degrees
$d_1,\ldots,d_{2m}$. If $x_{i_2}\neq x_{p_3},x_{p_4}$ and
$x_{i_3}\neq x_{p_3},x_{p_4}$, let
$G_1=G_0-\{u_2x_{i_2},u_2x_{i_3},u_1x_{p_3},u_1x_{p_4}\}+\{u_2x_{p_3},u_2x_{p_4},u_1x_{i_2},u_1x_{i_3}\}$,
then $G_1$ is a realization of $\pi$ and contains
$C_3,\ldots,C_{2m}$ on those vertices with degrees
$d_1,\ldots,d_{2m}$.

Now, we assume that $d_{2m+1}=2$ and $d_1=m+1$, or assume that
$d_{2m+1}=1$ and $d_1=m$.

Then $\pi=(d_1,m^{m+1},m-1,m-2,\ldots,2,d_{2m+1})$. Let
$\omega_1=(m^{m+1},(m-1)^2,m-2,\ldots,2)$,
$\omega_2=(m^2,(m-1)^m,m-2,\ldots,2)$ and
$\omega_3=((m-1)^{m+2},m-2,\ldots,3)$. Clearly, $\omega_1$ is the
residual sequence obtained from $\pi$ by laying off the term
$d_{2m+1}$, $\omega_2$ is the residual sequence obtained from
$\omega_1$ by laying off the term $m-1$ and $\omega_3$ is the
residual sequence obtained from $\omega_2$ by laying off the term
$2$. By Theorem 2.1 and Lemma 2.1, $\omega_3$ has a pancyclic
realization $G_3$. For $i=2,1,0$ in turn, let $G_i$ be the graph
obtained from $G_{i+1}$ by adding a new vertex $u_i$ so that $u_i$
is adjacent to those vertices whose degrees are reduced by one unity
in going from $\omega_i$ to $\omega_{i+1}$. Then $G_0$ is a
realization of $\pi$ and contains $G_3$. Let $C=x_1x_2\ldots
x_{2m-2}x_1$ be a Hamilton cycle of $G_3$. Then
$V(G_0)=\{x_1,\ldots,x_{2m-2},u_0,u_1,u_2\}$,
$d_{G_0}(u_0)=d_{2m+1}$, $d_{G_0}(u_1)=d_{m+3}=m-1$ and
$d_{G_0}(u_2)=d_{2m}=2$. Thus $C$ lies on those vertices with
degrees $d_1,\ldots,d_{m+2},d_{m+4},\ldots,d_{2m-1}$. For
convenience, let $\{i_1,\ldots,i_{m+1}\}\subseteq
\{1,2,\ldots,2m-2\}$ so that $d_{G_0}(x_{i_j})=d_j$ for $1\leq j\leq
m+1$. Clearly, $u_2x_{i_1},u_2x_{i_2}\in E(G_0)$ and $u_1x_{i_j}\in
E(G_0)$ for $3\leq j\leq m+1$. Moreover, there are four distinct
positive integers $p_1,p_2,p_3,p_4\in \{i_1,\ldots,i_{m+1}\}$ so
that $x_{p_1}x_{p_2},x_{p_3}x_{p_4}\in E(C)$. Without loss of
generality, if $x_{i_1}=x_{p_1}$ and $x_{i_2}=x_{p_2}$, then
$u_2x_{p_1}, u_2x_{p_2}, u_1x_{p_3}, u_1x_{p_4}\in E(G_0)$, and
hence $G_0$ contains $C_3,\ldots, C_{2m}$ on the vertices
$x_1,\ldots,x_{2m-2},u_1,u_2$. In other words, $G_0$ contains
$C_3,\ldots, C_{2m}$ on those vertices with degrees
$d_1,\ldots,d_{2m}$. If $x_{i_1}=x_{p_1}$ and $x_{i_2}\neq x_{p_2}$,
let $G_1=G_0-\{u_2x_{i_2},u_1x_{p_2}\}+\{u_2x_{p_2},u_1x_{i_2}\}$,
then $G_1$ is a realization of $\pi$ and
$u_2x_{p_1},u_2x_{p_2},u_1x_{p_3},u_1x_{p_4}\in E(G_1)$, implying
that $G_1$ contains $C_3,\ldots, C_{2m}$ on those vertices with
degrees $d_1,\ldots,d_{2m}$. If $x_{i_1}\not=x_{p_1},x_{p_2}$ and
$x_{i_2}\not=x_{p_1},x_{p_2}$, let
$G_1=G_0-\{u_2x_{i_1},u_2x_{i_2},u_1x_{p_1},u_1x_{p_2}\}+\{u_2x_{p_1},u_2x_{p_2},u_1x_{i_1},u_1x_{i_2}\}$,
then $G_1$ is a realization of $\pi$ and
$u_2x_{p_1},u_2x_{p_2},u_1x_{p_3},u_1x_{p_4}\in E(G_1)$. Thus $G_1$
contains $C_3,\ldots, C_{2m}$ on those vertices with degrees
$d_1,\ldots,d_{2m}$.\ \ $\Box$

{\bf Lemma 2.7}\ \ {\it Let $\ell\in \{7,8\}$, $n\geq\ell$ and
$\pi=(d_1,\ldots,d_n)\in GS_n$ with $d_{\ell-2}\geq 4$,
$d_{\ell-1}\geq 3$ and $d_\ell\geq 2$. Then $\pi$ has a realization
containing $C_3,\ldots,C_\ell$ on those vertices with degrees
$d_1,\ldots,d_\ell$.}

{\bf Proof.}\ \ Firstly, we consider the case that $\ell=8$. It
follows from Lemma 2.1 that Lemma 2.7 holds for $n=8$. Assume $n\geq
9$. If $d_8\geq 4$, by Lemma 2.2, then $\pi$ has a realization
containing $C_3,\ldots,C_8$ on those vertices with degrees
$d_1,\ldots,d_8$. Assume $d_8\leq 3$. If $d_{d_n}\geq 5$, then the
residual sequence $\pi_n'=(d_1',\ldots,d_{n-1}')$ obtained from
$\pi$ by laying off $d_n$ satisfies $d_6'\geq 4$, $d_7'\geq 3$ and
$d_8'\geq 2$. By the induction hypothesis, $\pi_n'$ has a
realization containing $C_3,\ldots,C_8$ on those vertices with
degrees $d_1',\ldots,d_8'$, and so $\pi$ has a realization
containing $C_3,\ldots,C_8$ on those vertices with degrees
$d_1,\ldots,d_8$. Assume $d_{d_n}=\cdots=d_6=4$ and $d_n\ge 1$.
Denote $p=\max\{i|d_i\geq 5\}$.

{\bf Case 1}\ \ $d_n=1$.

Then $p=0$. If $n\geq 10$, let
$\pi'=(d_1,\ldots,d_{n-2},d_{n-1}-1)$, if $n=9$ and $d_8\geq 3$, let
$\pi'=(d_1,\ldots,d_7,d_8-1)$, and if $n=9$, $d_8=2$ and $d_7=4$,
let $\pi'=(d_1,\ldots,d_6,d_7-1,d_8)$. It is easy to see that
$\pi'\in GS_{n-1}$. Denote $\pi'=(d_1',\ldots,d_{n-1}')$. Clearly,
$\pi'$ satisfies $d_6\geq 4$, $d_7\geq 3$ and $d_8\geq 2$. By the
induction hypothesis, $\pi'$ has a realization containing
$C_3,\ldots,C_8$ on those vertices with degrees $d_1',\ldots,d_8'$,
and so $\pi$ has a realization containing $C_3,\ldots,C_8$ on those
vertices with degrees $d_1,\ldots,d_8$. Thus we may assume that
$n=9$, $d_8=2$ and $d_7=3$, that is, $\pi=(4^6,3,2,1)$. Let
$\pi'=(4^6,3^2)$ be obtained from $\pi$ by deleting $d_9=1$ and
decreasing $d_8$ by one unity. Fig. 3 shows that $\pi'$ has a
pancyclic realization, implying that $\pi$ has a realization
containing $C_3,\ldots,C_8$ on those vertices with degrees
$d_1,\ldots,d_8$.

{\bf Case 2} \ \ $d_n=2$.

Then $p\le 1$. If $p=1$ and $n\geq 10$, or if $p=1$, $n=9$ and
$d_8\geq 3$, let $\pi'=(d_1-1,d_2,\ldots,d_{n-2},d_{n-1}-1)$. If
$p=1$, $n=9$, $d_8=2$ and $d_7=4$, let
$\pi'=(d_1-1,d_2,\ldots,d_6,d_7-1,d_8)$. If $p=0$ and $n\geq 11$, or
if $p=0$, $n=10$ and $d_8\geq 3$, let
$\pi'=(d_1,\ldots,d_{n-3},d_{n-2}-1,d_{n-1}-1)$. If $p=0$, $n=10$,
$d_8=2$ and $d_7$=4, let $\pi'=(d_1,\ldots,d_6,d_7-1,d_8,d_9-1)$. It
is easy to see that $\pi'\in GS_{n-1}$. Denote
$\pi'=(d_1',\ldots,d_{n-1}')$. Clearly, $\pi'$ satisfies $d_6\geq
4$, $d_7\geq 3$ and $d_8\geq 2$. By the induction hypothesis, $\pi'$
has a realization containing $C_3,\ldots,C_8$ on those vertices with
degrees $d_1',\ldots,d_8'$, and so $\pi$ has a realization
containing $C_3,\ldots,C_8$ on those vertices with degrees
$d_1,\ldots,d_8$. Thus we may assume that $p=1$, $n=9$, $d_8=2$ and
$d_7=3$, or assume that $p=0$, $n=10$, $d_8=2$ and $d_7=3$, or
assume that $p=0$ and $n=9$. If $p=0$, $n=10$, $d_8=2$ and $d_7=3$,
then $\pi=(4^6,3,2^3)$, which is impossible as $\sg(\pi)$ is odd. If
$p=1$, $n=9$, $d_8=2$ and $d_7=3$, then $\pi=(5,4^5,3,2^2)$ or
$(7,4^5,3,2^2)$. If $\pi=(5,4^5,3,2^2)$ (respectively,
$\pi=(7,4^5,3,2^2)$), let $\pi'=(4^6,2^2)$ (respectively,
$\pi'=(6,4^5,2^2)$) be obtained from $\pi$ by deleting $d_9=2$ and
decreasing $d_1$ and $d_7$ each by one unity. Fig. 3 (respectively,
Fig. 4) shows that $\pi'$ has a pancyclic realization. Thus $\pi$
has a realization containing $C_3,\ldots,C_8$ on those vertices with
degrees $d_1,\ldots,d_8$. If $p=0$ and $n=9$, then $\pi=(4^7,2^2)$
or $(4^6,3^2,2)$. If $\pi=(4^7,2^2)$ (respectively,
$\pi=(4^6,3^2,2)$), let $\pi'=(4^5,3^2,2)$ (respectively,
$\pi'=(4^6,2^2)$) be obtained from $\pi$ by deleting $d_9=2$ and
decreasing $d_6$ and $d_7$ each by one unity (respectively,
decreasing $d_7$ and $d_8$ each by one unity). Fig. 2 (respectively,
Fig. 3) shows that $\pi'$ has a pancyclic realization. Thus $\pi$
has a realization containing $C_3,\ldots,C_8$ on those vertices with
degrees $d_1,\ldots,d_8$.

{\bf Case 3}\ \ $d_n=3$.

Then $p\le 2$. If $p=2$, let
$\pi'=(d_1-1,d_2-1,d_3,\ldots,d_{n-2},d_{n-1}-1)$. If $p=1$ and
$n\geq 10$, or if $p=1$, $n=9$ and $d_7=4$, let
$\pi'=(d_1-1,d_2,\ldots,d_{n-3},d_{n-2}-1,d_{n-1}-1)$. If $p=0$ and
$n\geq 11$, or if $p=0$, $n=10$ and $d_7=4$, let
$\pi'=(d_1,\ldots,d_{n-4},d_{n-3}-1,d_{n-2}-1,d_{n-1}-1)$. It is
easy to see that $\pi'\in GS_{n-1}$. Denote
$\pi'=(d_1',\ldots,d_{n-1}')$. Clearly, $\pi'$ satisfies $d_6\geq
4$, $d_7\geq 3$ and $d_8\geq 2$. By the induction hypothesis, $\pi'$
has a realization containing $C_3,\ldots,C_8$ on those vertices with
degrees $d_1',\ldots,d_8'$, and so $\pi$ has a realization
containing $C_3,\ldots,C_8$ on those vertices with degrees
$d_1,\ldots,d_8$. Thus we may assume that $p=1$, $n=9$ and $d_7=3$,
or assume that $p=0$, $n=10$ and $d_7=3$, or assume that $p=0$ and
$n=9$. If $p=1$, $n=9$ and $d_7=3$, then $\pi=(5,4^5,3^3)$ or
$(7,4^5,3^3)$. If $\pi=(5,4^5,3^3)$ (respectively, $(7,4^5,3^3)$),
let $\pi'=(4^6,2^2)$ (respectively, $\pi'=(6,4^5,2^2)$) be obtained
from $\pi$ by deleting $d_9=3$ and decreasing $d_1$, $d_7$ and $d_8$
each by one unity. Fig. 3 (respectively, Fig. 4) shows that $\pi'$
has a pancyclic realization. Thus $\pi$ has a realization containing
$C_3,\ldots,C_8$ on those vertices with degrees $d_1,\ldots,d_8$. If
$p=0$, $n=10$ and $d_7=3$, then $\pi=(4^6,3^4)$. Let
$\pi'=(d_1',\ldots,d_9')=(4^6,2^3)$ be obtained from $\pi$ by
deleting $d_{10}=3$ and decreasing $d_7$, $d_8$ and $d_9$ each by
one unity, and let $\pi''=(4^4,3^2,2^2)$ be obtained from $\pi'$ by
deleting $d_9'=2$ and decreasing $d_5'$ and $d_6'$ each by one
unity. Fig. 1 shows that $\pi''$ has a pancyclic realization,
implying that $\pi$ has a realization containing $C_3,\ldots,C_8$ on
those vertices with degrees $d_1,\ldots,d_8$. If $p=0$ and $n=9$,
then $\pi=(4^7,3^2)$. Let $\pi'=(4^5,3^2,2)$ be obtained from $\pi$
by deleting $d_9=3$ and decreasing $d_6$, $d_7$ and $d_8$ each by
one unity. Fig. 2 shows that $\pi'$ has a pancyclic realization,
implying that $\pi$ has a realization containing $C_3,\ldots,C_8$ on
those vertices with degrees $d_1,\ldots,d_8$.

The proof of the case that $\ell=7$ is similar, we omit it here.\ \
$\Box$

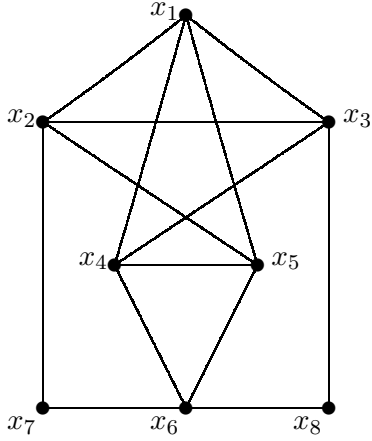
\begin{center}
\setlength{\unitlength}{0.95mm}
\begin{picture}(160,70)
\put(40,70){\circle*{2}}
 \put(35,70){$x_5$}
 \put(30,50){\circle*{2}}
 \put(25,50){$x_1$}
 \put(50,50){\circle*{2}}
 \put(52,50){$x_2$}
 \put(20,40){\circle*{2}}
 \put(15,40){$x_7$}
 \put(60,40){\circle*{2}}
 \put(62,40){$x_8$}
 \put(30,30){\circle*{2}}
 \put(25,28){$x_3$}
 \put(50,30){\circle*{2}}
 \put(52,28){$x_4$}
 \put(40,10){\circle*{2}}
 \put(35,10){$x_6$}
 \put(40,70){\line(1,-2){10}}
 \put(40,70){\line(-1,-2){10}}
 \put(20,40){\line(1,1){10}}
 \put(20,40){\line(1,-1){10}}
 \put(60,40){\line(-1,1){10}}
 \put(60,40){\line(-1,-1){10}}
 \put(30,50){\line(1,0){20}}
 \put(30,30){\line(1,0){20}}
 \put(40,10){\line(-1,2){10}}
 \put(40,10){\line(1,2){10}}
 \put(30,30){\line(1,1){20}}
 \put(50,30){\line(-1,1){20}}
 \qbezier(40,70)(100,40)(40,10)
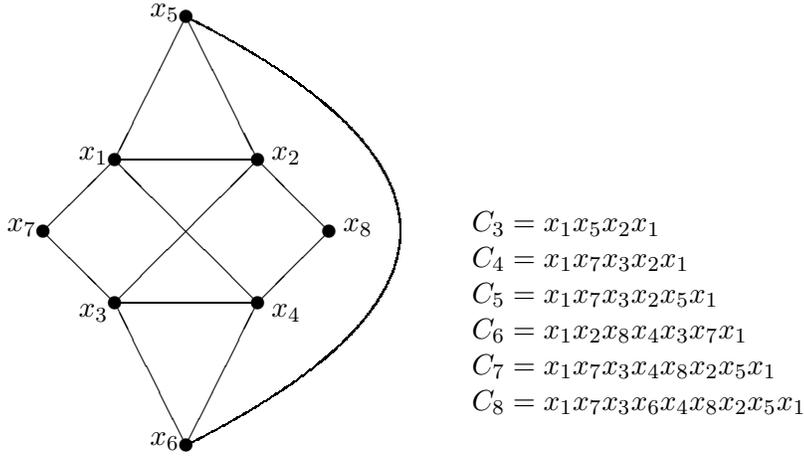
 \put(15,4){Fig. 1.\ \ A pancyclic realization of $\pi=(4^4,3^2,2^2)$.}
 \put(80,40){$C_3=x_1x_5x_2x_1$}
 \put(80,35){$C_4=x_1x_7x_3x_2x_1$}
 \put(80,30){$C_5=x_1x_7x_3x_2x_5x_1$}
 \put(80,25){$C_6=x_1x_2x_8x_4x_3x_7x_1$}
 \put(80,20){$C_7=x_1x_7x_3x_4x_8x_2x_5x_1$}
 \put(80,15){$C_8=x_1x_7x_3x_6x_4x_8x_2x_5x_1$}
\end{picture}
\end{center}

\begin{center}
\setlength{\unitlength}{0.95mm}
\begin{picture}(160,70)
\put(40,70){\circle*{2}}
 \put(35,70){$x_1$}
 \put(20,55){\circle*{2}}
 \put(15,55){$x_2$}
 \put(60,55){\circle*{2}}
 \put(62,55){$x_3$}
 \put(30,35){\circle*{2}}
 \put(25,35){$x_4$}
 \put(50,35){\circle*{2}}
 \put(52,35){$x_5$}
 \put(20,15){\circle*{2}}
 \put(15,12){$x_6$}
 \put(40,15){\circle*{2}}
 \put(35,12){$x_7$}
 \put(60,15){\circle*{2}}
 \put(55,12){$x_8$}
  \qbezier(40,70)(30,62)(20,55)
  \qbezier(40,70)(50,62)(60,55)
  \qbezier(40,70)(35,52)(30,35)
  \qbezier(40,70)(45,52)(50,35)
  \qbezier(20,55)(40,55)(60,55)
  \qbezier(20,55)(35,45)(50,35)
  \qbezier(60,55)(45,45)(30,35)
  \qbezier(20,55)(20,35)(20,15)
  \qbezier(60,55)(60,35)(60,15)
  \qbezier(30,35)(25,25)(20,15)
  \qbezier(50,35)(45,25)(40,15)
  \qbezier(20,15)(30,15)(40,15)
  \qbezier(40,15)(50,15)(60,15)
  \qbezier(30,35)(40,35)(50,35)
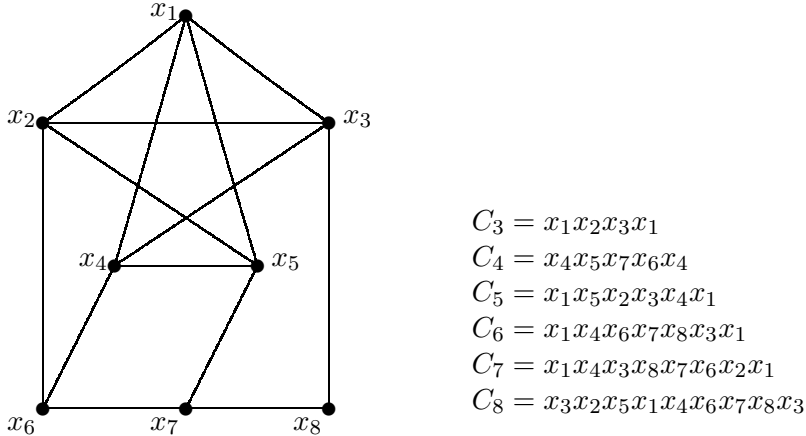
 \put(15,4){Fig. 2.\ \ A pancyclic realization of $\pi=(4^5,3^2,2)$.}
 \put(80,40){$C_3=x_1x_2x_3x_1$}
 \put(80,35){$C_4=x_4x_5x_7x_6x_4$}
 \put(80,30){$C_5=x_1x_5x_2x_3x_4x_1$}
 \put(80,25){$C_6=x_1x_4x_6x_7x_8x_3x_1$}
 \put(80,20){$C_7=x_1x_4x_3x_8x_7x_6x_2x_1$}
 \put(80,15){$C_8=x_3x_2x_5x_1x_4x_6x_7x_8x_3$}
\end{picture}
\end{center}

\begin{center}
\setlength{\unitlength}{0.95mm}
\begin{picture}(160,70)
\put(40,70){\circle*{2}}
 \put(35,70){$x_1$}
 \put(20,55){\circle*{2}}
 \put(15,55){$x_2$}
 \put(60,55){\circle*{2}}
 \put(62,55){$x_3$}
 \put(30,35){\circle*{2}}
 \put(25,35){$x_4$}
 \put(50,35){\circle*{2}}
 \put(52,35){$x_5$}
 \put(20,15){\circle*{2}}
 \put(15,12){$x_7$}
 \put(40,15){\circle*{2}}
 \put(35,12){$x_6$}
 \put(60,15){\circle*{2}}
 \put(55,12){$x_8$}
  \qbezier(40,70)(30,62)(20,55)
  \qbezier(40,70)(50,62)(60,55)
  \qbezier(40,70)(35,52)(30,35)
  \qbezier(40,70)(45,52)(50,35)
  \qbezier(20,55)(40,55)(60,55)
  \qbezier(20,55)(35,45)(50,35)
  \qbezier(60,55)(45,45)(30,35)
  \qbezier(20,55)(20,35)(20,15)
  \qbezier(60,55)(60,35)(60,15)
  \qbezier(30,35)(35,25)(40,15)
  \qbezier(50,35)(45,25)(40,15)
  \qbezier(20,15)(30,15)(40,15)
  \qbezier(40,15)(50,15)(60,15)
  \qbezier(30,35)(40,35)(50,35)
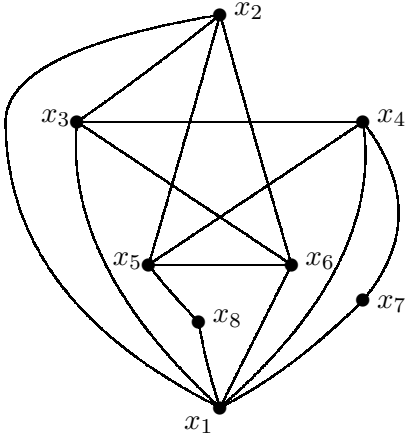
 \put(15,4){Fig. 3.\ \ A pancyclic realization of $\pi=(4^6,2^2)$.}
 \put(80,40){$C_3=x_4x_5x_6x_4$}
 \put(80,35){$C_4=x_1x_4x_6x_5x_1$}
 \put(80,30){$C_5=x_1x_4x_3x_2x_5x_1$}
 \put(80,25){$C_6=x_1x_2x_7x_6x_8x_3x_1$}
 \put(80,20){$C_7=x_4x_5x_2x_7x_6x_8x_3x_4$}
 \put(80,15){$C_8=x_2x_5x_1x_4x_3x_8x_6x_7x_2$}
\end{picture}
\end{center}

\begin{center}
\setlength{\unitlength}{0.95mm}
\begin{picture}(160,70)
\put(40,70){\circle*{2}}
 \put(42,70){$x_2$}
 \put(20,55){\circle*{2}}
 \put(15,55){$x_3$}
 \put(60,55){\circle*{2}}
 \put(62,55){$x_4$}
 \put(30,35){\circle*{2}}
 \put(25,35){$x_5$}
 \put(50,35){\circle*{2}}
 \put(52,35){$x_6$}
 \put(40,15){\circle*{2}}
 \put(35,12){$x_1$}
 \put(60,30){\circle*{2}}
 \put(62,29){$x_7$}
 \put(37,27){\circle*{2}}
 \put(39,27){$x_8$}
  \qbezier(40,70)(30,62)(20,55)
  \qbezier(40,70)(35,52)(30,35)
  \qbezier(40,70)(45,52)(50,35)
  \qbezier(20,55)(40,55)(60,55)
  \qbezier(20,55)(35,45)(50,35)
  \qbezier(60,55)(45,45)(30,35)
  \qbezier(50,35)(45,25)(40,15)
  \qbezier(30,35)(40,35)(50,35)
  \qbezier(40,15)(63,35)(60,55)
  \qbezier(40,15)(18,35)(20,55)
  \qbezier(40,15)(38,21)(37,27)
  \qbezier(30,35)(33,31)(37,27)
  \qbezier(40,15)(50,20)(60,30)
  \qbezier(60,30)(70,42)(60,55)
  \qbezier(40,15)(10,30)(10,55)
  \qbezier(10,55)(10,65)(40,70)
 \put(15,4){Fig. 4.\ \ A pancyclic realization of $\pi=(6,4^5,2^2)$.}
 \put(80,40){$C_3=x_1x_2x_3x_1$}
 \put(80,35){$C_4=x_1x_8x_5x_6x_1$}
 \put(80,30){$C_5=x_2x_5x_4x_3x_6x_2$}
 \put(80,25){$C_6=x_4x_3x_6x_5x_8x_1x_4$}
 \put(80,20){$C_7=x_4x_3x_6x_5x_8x_1x_7x_4$}
 \put(80,15){$C_8=x_5x_2x_6x_3x_4x_7x_1x_8x_5$}
\end{picture}
\end{center}

{\bf Lemma 2.8}\ \ {\it Let $m\geq 4$, $n\geq 2m$ and
$\pi=(d_1,\ldots,d_n)\in GS_n$ with $d_{2m+1-i}\geq i+1$ for $1\leq
i\leq m-1$. Then $\pi$ has a realization containing $C_3,\ldots,
C_{2m}$ on those vertices with degrees $d_1,\ldots,d_{2m}$.}

{\bf Proof.}\ \ It follows from Lemma 2.7 that Lemma 2.8 holds for
$m=4$. We assume $m\ge 5$. It follows from Lemma 2.1 that Lemma 2.8
holds for $n=2m$. We further assume $n\ge 2m+1$ and $d_n\ge 1$. If
$d_{2m}\geq m$, by Lemma 2.2, then $\pi$ has a realization
containing $C_3,\ldots, C_{2m}$ on those vertices with degrees
$d_1,\ldots,d_{2m}$. Assume $d_{2m}\leq m-1$.

If $d_{d_n} \geq m+1$, then the residual sequence
$\pi_n'=(d_1',\ldots,d_{n-1}')$ obtained from $\pi$ by laying off
$d_n$ satisfies $d_{2m+1-i}'\geq i+1$ for $1\leq i\leq m-1$. By the
induction hypothesis, $\pi_n'$ has a realization containing
$C_3,\ldots,C_{2m}$ on those vertices with degrees
$d_1',\ldots,d_{2m}'$, implying that $\pi$ has a realization
containing $C_3,\ldots,C_{2m}$ on those vertices with degrees
$d_1,\ldots,d_{2m}$. Thus we may assume that
$d_{d_n}=\cdots=d_{m+2}=m$. Denote $p=\max\{i|d_i\geq m+1\}$.

If $d_{d_n-1}\geq m+1$ and $n\geq 2m+2$, let
$\omega=(d_1-1,\ldots,d_{d_n-1}-1,d_{d_n},\ldots,d_{n-2},d_{n-1}-1)$,
and if $d_{d_n-1}\geq m+1$, $n=2m+1$ and there exists an integer $k$
with $1\le k\le m-2$ so that $d_{2m+1-k}\geq k+2$ (we choose $k$ to
be minimum), let
$\omega=(d_1-1,\ldots,d_{d_n-1}-1,d_{d_n},\ldots,d_{2m-k},d_{2m+1-k}-1,d_{2m+2-k},\ldots,d_{n-1})$.
By Lemma 2.3, $\omega\in GS_{n-1}$. Denote
$\omega=(\rho_1,\ldots,\rho_{n-1})$. Clearly, $\omega$ satisfies
$\rho_{2m+1-i}\geq i+1$ for $1\leq i\leq m-1$. By the induction
hypothesis, $\omega$ has a realization containing
$C_3,\ldots,C_{2m}$ on those vertices with degrees
$\rho_1,\ldots,\rho_{2m}$, and so $\pi$ has a realization containing
$C_3,\ldots,C_{2m}$ on those vertices with degrees
$d_1,\ldots,d_{2m}$. If $d_{d_n-1}\geq m+1$, $n=2m+1$ and
$d_{2m+1-i}=i+1$ for $1\le i\le m-2$, by $1\le d_{2m+1}\le 2$ and
Lemma 2.6, then $\pi$ has a realization containing $C_3,\ldots,
C_{2m}$ on those vertices with degrees $d_1,\ldots,d_{2m}$. Thus we
may further assume that $d_{d_n-1}=m$.

If $d_{m+3}=\cdots=d_{2m+1}=m-1$, let
$\rho=(d_1,\ldots,d_{m+3},d_{m+4}-1,\ldots,d_{2m}-1,d_{2m+2}-1,d_{2m+3}-1,d_{2m+4},\ldots,d_n)$
when $p=0$ and $n\geq 2m+3$,
$\rho=(d_1-1,d_2,\ldots,d_{m+3},d_{m+4}-1,\ldots,d_{2m}-1,d_{2m+2}-1,d_{2m+3},\ldots,d_n)$
when $p=1$ and $n\geq 2m+2$,
$\rho=(d_1-1,d_2-1,d_3,\ldots,d_{m+3},d_{m+4}-1,\ldots,d_{2m}-1,d_{2m+2},\ldots,d_n)$
when $p=2$, and
$\rho=(d_1-1,\ldots,d_p-1,d_{p+1},\ldots,d_{m+3},d_{m+4}-1,\ldots,d_{2m+2-p}-1,d_{2m+3-p},\ldots,d_{2m},d_{2m+2},\ldots,d_n)$
when $p\ge 3$, by Lemma 2.4, then $\rho$ is graphic. Let
$\rho'=(\rho_1,\ldots,\rho_{n-1})$ be a rearrangement in
non-increasing order of the $n-1$ terms of $\rho$. Clearly, $\rho'$
satisfies $\rho_{2m+1-i}\geq i+1$ for $1\leq i\leq m-1$. By the
induction hypothesis, $\rho'$ has a realization $G'$ containing
$C_3,\ldots,C_{2m}$ on those vertices with degrees
$\rho_1,\ldots,\rho_{2m}$. A realization $G$ of $\pi$ can be
obtained from $G'$ by adding a new vertex that is adjacent to those
vertices whose degrees are reduced by one unity in going from $\pi$
to $\rho'$. Let $H$ be the subgraph of $G'$ induced by those
vertices with degrees $\rho_1,\ldots,\rho_{2m}$. Clearly, $H$ is a
subgraph of $G$ and is pancyclic. By Theorem 2.3, $\pi$ has a
realization containing $H$ (and so containing $C_3,\ldots,C_{2m}$)
on those vertices with degrees $d_1,\ldots,d_{2m}$. Thus we now may
further assume that $d_{m+3}=m$, or assume that $d_{m+3}=m-1$ and
$d_{2m+1}\leq m-2$, or assume that $d_{m+3}=\cdots=d_{2m+1}=m-1$ and
$n\leq 2m+2-p$ for $p\in \{0,1\}$.

Clearly, $d_{2m}\ge p+2$. Let
$\omega_1=(d_1-1,\ldots,d_m-1,d_{m+1},d_{m+3},\ldots,d_{n})$. If
$p\ge 1$ or $d_{2m}\le m-2$, let
$\omega_2=(d_1-2,\ldots,d_p-2,d_{p+1}-1,\ldots,d_{m+1}-1,d_{m+3}-1,\ldots,d_r-1,d_{r+1},\ldots,d_{2m-1},d_{2m+1},\ldots,d_{n})$,
where $r=m+1+d_{2m}-p$. If $p=0$ and $d_{2m}=m-1$, let
$\omega_2=(d_1-2,\ldots,d_p-2,d_{p+1}-1,\ldots,d_{m+1}-1,d_{m+3}-1,\ldots,d_{2m-1}-1,d_{2m+1}-1,d_{2m+2},\ldots,d_{n})$.
By Lemma 2.5, $\omega_2$ is graphic. Let $\rho=(\rho_1,\ldots,
\rho_{n-2})$ be a rearrangement in non-increasing order of the $n-2$
terms of $\omega_2$. Clearly, $\rho_{2m-1-i}\geq i+1$ for $1\leq
i\leq m-2$ and
$\{d_1-2,\ldots,d_p-2,d_{p+1}-1,\ldots,d_{m+1}-1\}\subseteq
\{\rho_1,\ldots,\rho_{2m-2}\}$. Moreover, it is easy to see that
$d_{m+3}-1\in \{\rho_1,\ldots,\rho_{2m-2}\}$. By the induction
hypothesis, $\rho$ has a realization $G_2$ contains $C_3, \ldots ,
C_{2m-2}$ on those vertices with degrees $\rho_1, \ldots ,
\rho_{2m-2}$. For $i=1,0$ in turn, let $G_i$ be the graph obtained
from $G_{i+1}$ by adding a new vertex $u_i$ that is adjacent to
those vertices whose degree are reduced by one unity in going from
$\omega_{i}$ to $\omega_{i+1}$, where $\omega_0=\pi$. Then $G_0$ is
a realization of $\pi$ and contains $G_2$. Let $C=x_1x_2\cdots
x_{2m-2}x_1$ be the $C_{2m-2}$ in $G_2$. Then $C$ contains $m+2$
vertices, say $x_{i_1},\ldots,x_{i_{m+2}}$, so that
$d_{G_0}(x_{i_j})=d_j$ for $1\leq j\leq m+1$ and
$d_{G_0}(x_{i_{m+2}})=d_{m+3}$. Moreover, $d_{G_0}(u_0)=d_{m+2}=m$
and $d_{G_0}(u_1)=d_{2m}$. Clearly,
$u_1x_{i_{m+1}},u_1x_{i_{m+2}}\in E(G_0)$, $u_1x_{i_j}\in E(G_0)$
for $1\leq j\leq p$ and $u_0x_{i_j}\in E(G_0)$ for $1\leq j\leq m$.
There are four distinct positive integers $p_1,p_2,p_3,p_4\in
\{i_1,\ldots,i_{m+2}\}$ so that $x_{p1}x_{p2},x_{p3}x_{p4}\in E(C)$.
We now show that $\pi$ has a realization $G$ containing
$C_3,\ldots,C_{2m}$ on the vertices $x_1,\ldots,x_{2m-2},u_0,u_1$.
Without loss of generality, we only need to consider three cases.

{\bf Case 1}\ \ $x_{p_1},x_{p_2}\in
\{x_{i_{p+1}},\ldots,x_{i_{m+2}}\}$.

If $x_{p_1}=x_{i_{m+1}}$ and $x_{p_2}=x_{i_{m+2}}$, that
$u_1x_{p_1},u_1x_{p_2},u_0x_{p_3},u_0x_{p_4}\in E(G_0)$, implying
that $G_0$ contains $C_3,\ldots,C_{2m}$ on the vertices
$x_1,\ldots,x_{2m-2},u_0,u_1$. If $x_{p_1}=x_{i_{m+1}}$ and
$x_{p_2}\neq x_{i_{m+2}}$, let
$G_1=G_0-\{u_0x_{i_{p_2}},u_1x_{i_{m+2}}\}+\{u_0x_{x_{i_{m+2}}},u_1x_{p_2}\}$,
then $G_1$ is a realization of $\pi$,
$u_1x_{p_1},u_1x_{p_2},u_0x_{p_3},u_0x_{p_4}\in E(G_1)$, $G_1$
contains $C_3,\ldots,C_{2m}$ on the vertices
$x_1,\ldots,x_{2m-2},u_0,u_1$. If $x_{p_1}\neq
x_{i_{m+1}},x_{i_{m+2}}$ and $x_{p_2}\neq x_{i_{m+1}},x_{i_{m+2}}$,
let
$$G_1=G_0-\{u_0x_{p_1},u_0x_{p_2},u_1x_{i_{m+1}},u_1x_{i_{m+2}}\}+\{u_0x_{i_{m+1}},u_0x_{i_{m+2}},u_1x_{p_1},u_1x_{p_2}\},$$
then $u_1x_{p_1},u_1x_{p_2},u_0x_{p_3},u_0x_{p_4}\in E(G_1)$, and
$G_1$ is a realization of $\pi$ and contains $C_3,\ldots,C_{2m}$ on
the vertices $x_1,\ldots,x_{2m-2},u_0,u_1$.

{\bf Case 2}\ \ $x_{p_1},x_{p_2}\in \{x_{i_1},\ldots,x_{i_p}\}$.

If $x_{p_3}\neq x_{i_{m+1}},x_{i_{m+2}}$ and $x_{p_4}\neq
x_{i_{m+1}},x_{i_{m+2}}$, then
$u_1x_{p_1},u_1x_{p_2},u_0x_{p_3},u_0x_{p_4}\in E(G_0)$ and $G_0$
contains $C_3,\ldots,C_{2m}$ on the vertices
$x_1,\ldots,x_{2m-2},u_0,u_1$. If $x_{p_3}\neq
x_{i_{m+1}},x_{i_{m+2}}$ and $x_{p_4}=x_{i_{m+1}}$, by $p\le m-3$,
let $G_1=G_0-\{u_0x_{i_m},u_1x_{p_4}\}+\{u_0x_{p_4},u_1x_{i_m}\}$,
then $G_1$ is a realization of $\pi$ and
$u_1x_{p_1},u_1x_{p_2},u_0x_{p_3},u_0x_{p_4}\in E(G_1)$. Thus $G_1$
contains $C_3,\ldots,C_{2m}$ on the vertices
$x_1,\ldots,x_{2m-2},u_0,u_1$. If $x_{p_3}=x_{i_{m+1}}$ and
$x_{p_4}=x_{i_{m+2}}$, then
$u_1x_{p_3},u_1x_{p_4},u_0x_{p_1},u_0x_{p_2}\in E(G_0)$ and $G_0$
contains $C_3,\ldots,C_{2m}$ on the vertices
$x_1,\ldots,x_{2m-2},u_0,u_1$.

{\bf Case 3}\ \ $x_{p_1}\in \{x_{i_1},\ldots,x_{i_p}\}$ and
$x_{p_2}\in \{x_{i_{p+1}},\ldots,x_{i_{m+2}}\}$.

If $x_{p_2}=x_{i_{m+1}}$ and $x_{p_3},x_{p_4}\not=x_{i_{m+2}}$, then
$u_1x_{p_1},u_1x_{p_2},u_0x_{p_3},u_0x_{p_4}\in E(G_0)$ and $G_0$
contains $C_3,\ldots,C_{2m}$ on the vertices
$x_1,\ldots,x_{2m-2},u_0,u_1$. If $x_{p_2}=x_{i_{m+1}}$ and
$x_{p_3}=x_{i_{m+2}}$, let
$G_1=G_0-\{u_0x_{i_m},u_1x_{p_3}\}+\{u_0x_{p_3},u_1x_{i_m}\}$ when
$x_{p_4}\not=x_{i_m}$ and
$G_1=G_0-\{u_0x_{i_{m-1}},u_1x_{p_3}\}+\{u_0x_{p_3},u_1x_{i_{m-1}}\}$
when $x_{p_4}=x_{i_m}$, then
$u_1x_{p_1},u_1x_{p_2},u_0x_{p_3},u_0x_{p_4}\in E(G_1)$ and $G_1$ is
a realization of $\pi$ and contains $C_3,\ldots,C_{2m}$ on the
vertices $x_1,\ldots,x_{2m-2},u_0,u_1$. If $x_{p_2}\neq
x_{i_{m+1}},x_{i_{m+2}}$, $x_{p_3}=x_{i_{m+1}}$ and
$x_{p_4}=x_{i_{m+2}}$, then
$u_0x_{p_1},u_0x_{p_2},u_1x_{p_3},u_1x_{p_4}\in E(G_0)$ and $G_0$
contains $C_3,\ldots,C_{2m}$ on the vertices
$x_1,\ldots,x_{2m-2},u_0,u_1$. If $x_{p_2}\neq
x_{i_{m+1}},x_{i_{m+2}}$, $x_{p_3}=x_{i_{m+1}}$ and $x_{p_4}\neq
x_{i_{m+2}}$, let
$G_1=G_0-\{u_1x_{p_3},u_0x_{p_2}\}+\{u_1x_{p_2},u_0x_{p_3}\}$, then
$u_1x_{p_1},u_1x_{p_2},u_0x_{p_3},u_0x_{p_4}\in E(G_1)$ and $G_1$ is
a realization of $\pi$ and contains $C_3,\ldots,C_{2m}$ on the
vertices $x_1,\ldots,x_{2m-2},u_0,u_1$. If $x_{p_2}\neq
x_{i_{m+1}},x_{i_{m+2}}$, $x_{p_3}\neq x_{i_{m+1}},x_{i_{m+2}}$ and
$x_{p_4}\neq x_{i_{m+1}},x_{i_{m+2}}$, let
$G_1=G_0-\{u_1x_{i_{m+1}},u_0x_{p_2}\}+\{u_1x_{p_2},u_0x_{i_{m+1}}\}$,
then $u_1x_{p_1},u_1x_{p_2},u_0x_{p_3},u_0x_{p_4}\in E(G_1)$ and
$G_1$ is a realization of $\pi$ and contains $C_3,\ldots,C_{2m}$ on
the vertices $x_1,\ldots,x_{2m-2},u_0,u_1$.

Let $H$ be the subgraph of $G$ induced by
$\{x_1,\ldots,x_{2m-2},u_0,u_1\}$. Clearly, $H$ is pancyclic. By
Theorem 2.3, $\pi$ has a realization containing $H$ (and so
containing $C_3,\ldots,C_{2m}$) on those vertices with degrees
$d_1,\ldots,d_{2m}$. \ \ $\Box$

{\bf Lemma 2.9}\ \ {\it Let $m\geq 4$, $n\geq 2m+1$ and
$\pi=(d_1,\ldots,d_n)\in GS_n$ with $d_{2m+2-i}\geq i+1$ for $1\leq
i\leq m$. Then 
$\pi$ is potentially $_3C_{2m+1}$-graphic.}

{\bf Proof.}\ \ It follows from Lemma 2.1 that Lemma 2.9 holds for
$n=2m+1$. Assume $n\geq 2m+2$. If $d_{2m+1}\geq m+1$, by Theorem
1.2, then $\pi$ is potentially $_3C_{2m+1}$-graphic. Assume
$d_{2m+1}\le m$. If $d_{d_n} \geq m+2$, then the residual sequence
$\pi_n'=(d_1',\ldots,d_{n-1}')$ obtained from $\pi$ by laying off
$d_n$ satisfies $d_{2m+2-i}'\geq i+1$ for $1\leq i\leq m$. By the
induction hypothesis, $\pi_n'$ is potentially $_3C_{2m+1}$-graphic,
implying that $\pi$ is potentially $_3C_{2m+1}$-graphic. Assume
$d_{d_n}=\cdots=d_{m+2}=m+1$. Clearly, the residual sequence
$\pi_{m+2}'=(d_1',\ldots,d_{n-1}')$ obtained from $\pi$ by laying
off $d_{m+2}=m+1$ satisfies $d_{2m+1-i}'\geq i+1$ for $1\leq i\leq
m-1$. By Lemma 2.8, $\pi_{m+2}'$ has a realization $G'$ containing
$C_3,\ldots,C_{2m}$ on those vertices with degree
$d_1',\ldots,d_{2m}'$. A realization $G$ of $\pi$ can be obtained
from $G'$ by adding a new vertex $u$ that is adjacent to the
vertices whose degrees are reduced by one unity in going from $\pi$
to $\pi_{m+2}'$. By $d_{2m+1}\le m$, we can see that $C_{2m}$ lies
on those vertices with degrees
$d_1,\ldots,d_{m+1},d_{m+3},\ldots,d_{2m+1}$ and $u$ is adjacent to
$m+1$ vertices of $C_{2m}$. So $G$ contains $C_{2m+1}$. Thus $\pi$
is potentially $_3C_{2m+1}$-graphic.\ \ $\Box$

{\bf Proof of Theorem 1.3.}\ \ Theorem 1.3 follows immediately from
Theorem 1.1 and Lemmas 2.7--2.9.\ \ $\Box$

As an application of Theorem 1.3, we now present a proof of Theorem
1.4.

{\bf Proof of Theorem 1.4.}\ \ If $\ell\ge 5$, let
$\pi=(n-1,(\ell-2)^{\ell-2}, 1^{n-\ell+1})$, then $\pi$ is not
potentially $C_\ell$-graphic, and so $\sg(C_\ell,n)\ge
\sg(\pi)+2=2n+\ell^2-5\ell+6$. Moreover, for $m\ge 2$, since
$\pi=((n-1)^m,m^{n-m})$ (respectively,
$\pi=((n-1)^m,(m+1)^2,m^{n-m-2})$) is not potentially $C_{2m+1}$
(respectively, $C_{2m+2}$)-graphic, we have that $\sg(C_{2m+1},n)\ge
\sg(\pi)+2=m(2n-m-1)+2$ (respectively, $\sg(C_{2m+2},n)\ge
\sg(\pi)+2=m(2n-m-1)+4$). Therefore, $\sg(C_{2m+1},n)\ge
\max\{2n+(2m+1)^2-5(2m+1)+6,m(2n-m-1)+2\}=\max\{2n+4m^2-6m,m(2n-m-1)\}+2$
and $\sg(C_{2m+2},n)\ge
\max\{2n+(2m+2)^2-5(2m+2)+6,m(2n-m-1)+4\}=\max\{2n+4m^2-2m-2,m(2n-m-1)+2\}+2$.

In order to prove that $\sg(C_{2m+1},n)\le
\max\{2n+4m^2-6m,m(2n-m-1)\}+2$ for $n\ge 2m+1$, we only need to
show that if $n\ge 2m+1$ and $\pi=(d_1,\ldots,d_n)\in GS_n$ with
$\sg(\pi)\ge \max\{2n+4m^2-6m,m(2n-m-1)\}+2$, then $\pi$ is
potentially $C_{2m+1}$-graphic. If $d_{2m+2-i}\geq i+1$ for $1\leq
i\leq m$, by Theorem 1.3, then $\pi$ is potentially
$C_{2m+1}$-graphic. Assume that there is an $k$ with $1\leq k\leq m$
so that $d_{2m+2-k}\leq k$. Denote $f(k)=3k^2+(2n-8m-3)k+4m^2+2m$.
Then $\sg(\pi)=(d_1+\cdots+d_{2m+1-k})+d_{2m+2-k}+\cdots+d_n\le
(2m+1-k)(2m-k)+\sum\limits_{i=2m+2-k}^{n}\min\{2m+1-k,d_i\}+d_{2m+2-k}+\cdots+d_n=(2m+1-k)(2m-k)+2(d_{2m+2-k}+\cdots+d_n)\le
(2m+1-k)(2m-k)+2k(n-2m-1+k)=f(k)\le
\max\{f(1),f(m)\}=\max\{2n+4m^2-6m,m(2n-m-1)\}<\sg(\pi)$, a
contradiction.

Similarly, to show that $\sg(C_{2m+2},n)\le
\max\{2n+4m^2-2m-2,m(2n-m-1)+2\}+2$ for $n\ge 2m+2$, we let $n\ge
2m+2$ and $\pi=(d_1,\ldots,d_n)\in GS_n$ with $\sg(\pi)\ge
\max\{2n+4m^2-2m-2,m(2n-m-1)+2\}+2$. If $d_{2m+3-i}\geq i+1$ for
$1\leq i\leq m$, by Theorem 1.3, then $\pi$ is potentially
$C_{2m+2}$-graphic. Assume that there is an $k$ with $1\leq k\leq m$
so that $d_{2m+3-k}\leq k$. Denote $g(k)=3k^2+(2n-8m-7)k+4m^2+6m+2$.
Then $\sg(\pi)=(d_1+\cdots+d_{2m+2-k})+d_{2m+3-k}+\cdots+d_n\le
(2m+2-k)(2m+1-k)+\sum\limits_{i=2m+3-k}^{n}\min\{2m+2-k,d_i\}+d_{2m+3-k}+\cdots+d_n=(2m+2-k)(2m+1-k)+2(d_{2m+3-k}+\cdots+d_n)\le
(2m+2-k)(2m+1-k)+2k(n-2m-2+k)=g(k)\le
\max\{g(1),g(m)\}=\max\{2n+4m^2-2m-2,m(2n-m-1)+2\}<\sg(\pi)$, a
contradiction.\ \ $\Box$

\vskip 0.1cm

\noindent{\bf References}

\bref{[1]}{J.A. Bondy and U.S.R. Murty,}{ Graph Theory with
Applications,}{} { The Macmillan Press, London, 1976.}

\vskip 0.1cm

\bref{[2]}{Gang Chen et al., On potentially $_3C_5$-graphic
sequences,} { J. Combin. Math. Combin. Comput.,}{ 61}{ (2007),
141--148.}

\vskip 0.1cm

\bref{[3]}{Gang Chen and Ying-Mei Fan, On potentially
$_3C_6$-graphic sequences,} { J. Guangxi Normal Univ. (Nat. Sci.
Ed.),}{ 24(3)}{ (2006), 26--29.}

\bref{[4]} {P. Erd\H{o}s and T. Gallai, Graphs with prescribed
degrees of vertices (Hungarian),}{ Mat. Lapok,} { 11}{ (1960),
264--274.}

\vskip 0.1cm

\bref{[5]}{R.J. Gould, M.S. Jacobson and J. Lehel, Potentially
$G$-graphical degree sequences,}{ Combinatorics, Graph Theory, and
Algorithms,}{}{ Vol. I, II (Kalamazoo, MI, 1996), 451--460, New
Issues Press, Kalamazoo, MI, 1999.}

\vskip 0.1cm

\bref{[6]}{D.J. Kleitman and D.L. Wang, Algorithm for constructing
graphs and digraphs with given valences and factors,}{ Discrete
Math.,}{ 6}{ (1973), 79--88.}

\vskip 0.1cm

\bref{[7]}{Chunhui Lai, The smallest degree sum that yields
potentially $C_k$-graphical sequences,}{ J. Combin. Math. Combin.
Comput.,}{ 49}{ (2004), 57--64.}

\vskip 0.1cm

\bref{[8]}{J.S. Li and J.H. Yin, Extremal graph theory and degree
sequences,}{ Adv. Math.,}{ 33(3)}{ (2004), 273--283.}

\vskip 0.1cm

\bref{[9]}{C.St.J.A. Nash-Williams, Valency sequences which force
graphs to have hamiltonian circuits, Interim report, University of
Waterloo, Waterloo, 1970.}{}{}{}

\vskip 0.1cm

\bref{[10]}{E.F. Schmeichel and S.L. Hakimi, Pancyclic graphs and a
conjecture of Bondy and Chv\'{a}tal,} { J. Combin. Theory (B),}{
17}{ (1974), 22--34.}

\vskip 0.1cm

\bref{[11]}{J.H. Yin, A characterization for a graphic sequence to
be potentially $C_r$-graphic,}{ Sci. China Math.,}{ 53(11)}{ (2010),
2893--2905.}

\vskip 0.1cm

\bref{[12]}{J.H. Yin, K. Ye and J.Y. Li, Graphic sequences with a
realization containing cycles $C_3,\ldots,C_\ell$,}{ Appl. Math.
Comput.,}{ 353}{ (2019), 88--94.}

\end{document}